# Some Combinatorics behind Proofs

A. CARBONE [*]

November 28, 1995

**Abstract**

We try to bring to light some combinatorial structure underlying formal proofs in logic. We do this through the study of the Craig Interpolation Theorem which is properly a statement about the structure of formal derivations. We show that there is a generalization of the interpolation theorem to much more naive structures about sets, and then we show how both classical and intuitionistic versions of the statement follow by interpreting properly the set-theoretic language.

## 1 Introduction

There are interesting combinatorics behind formalized proofs in logic. In order to bring this out we need to find ways to extract combinatorial essence from proofs, to find 'homomorphisms' from proofs into combinatorics. We will describe here one way to do this, in connection with the Craig Interpolation Theorem.

Let us start with a combinatorial picture. Let $A$ and $B$ be two sets of points with a certain structure induced on them (think for instance of sets of points in an euclidean space with a certain distance function, or of a graph whose vertices are points of a set) and suppose you are given a finite number of operators that will allow you to define a pseudomap [1] $f$ between points of $A$ and $B$ (as discussed in section 2). This pseudomap does not have to be either surjective or everywhere defined. Namely there might be points in $A$ which might not be linked to any of the points in $B$ and viceversa. The theorem we present says essentially that there is a set $I$ of points (called the *interpolant* set) and two pseudomaps $f_A, f_B$ definable through the given set of operators which partially map $A$ surjectively into $I$ and $I$ into $B$. [2] Moreover, given any point $z \in I$ there are always points $x \in A$ and $y \in B$ such that $z \in f_A(x)$, $y \in f_B(z)$ and $y \in f(x)$.

The precise role of the operators mentioned above (from which the pseudomap $f$ is obtained) will be spelled out in the next section. For now let us say that the sets $A$ and $B$ and the pseudomap $f$ are obtained from some sets of 'trivial structure' by applying these operators a finite number of times, with the notion of trivial structure to be explained in section 2 also.

---

[*]Partially supported by the Lise-Meitner Stipendium # M00187-MAT (Austrian FWF.)

[1]A *pseudomap* $f$ is a relation between $A$ and $B$ in the sense of set theory. Since in everyday language the word *relation* has overtones of symmetry, we prefer to speak about pseudomaps. The shorthand $f(x)$ will refer to any point of $B$ that is $f$-related to $x \in A$.

[2]The sentence '$f$ partially maps $A$ into $I$' means that $f$ needs only be defined on a subset of $A$.

The theorem is essentially a geometrical formulation of the well-known Craig Interpolation Theorem which says that given two formulas $A$ and $B$ such that $A \rightarrow B$ is provable, there is a formula $I$ called an interpolant for $A$ and $B$ such that $I$ can be expressed in the language common to $A$ and $B$ – roughly speaking, all the symbols that appear in $I$ also appear in both $A$ and $B$ – and such that $A \rightarrow I$ and $I \rightarrow B$ are provable.

Our aim in this generalization is to bring out the mathematical structure lying below the proof of Craig's Theorem. The idea is to obtain Craig's Theorem by rereading points and sets in a logical language, namely by interpreting atomic occurrences in a proof as points in a set, formulas as sets and rules as operators on sets. The theorem is independent of the structure of the sets we consider. To establish the result only some compatibility of structures should be satisfied and this requirement will be captured by the properties of the operators.

As will be explained in section 5, this analysis was partly motivated by the close relation that exists between complexity bounds of interpolants and the problem of knowing whether $NP \cap CO-NP$ is contained in $P/poly$. We begin with some definitions in the next section which provide a combinatorial description of the logical notions of formulas, sequents, and rules of inference.

The author wishes to express her sincere thanks to Stephen Semmes for the uncountable number of comments which greatly improved the readability of this paper.

## 2  Spaces and Transformations

A *structured set* $S$ is a bipartite collection $\langle s_1 \ldots s_k | s_{k+1} \ldots s_n \rangle$ of sets of points where each set $s_i$ may be equipped with some additional structure and more than one copy of a set together with its structure might belong to $S$. We leave the reader free to think of his own favorite way to give a structure to a set of points. For example we can think of the sets $s_i$'s being finite graphs (as in figure (a)) or tessellations, i.e. geometric surfaces divided into nonoverlapping congruent polygons (as in (b)).

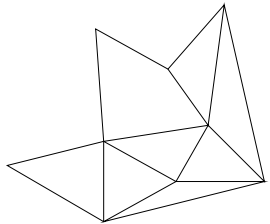
(a)

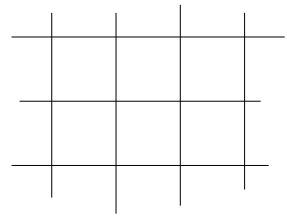
(b)

We say that the sets $s_1 \ldots s_k$ belong to the *first* component of $S$ and $s_{k+1} \ldots s_n$ to the *second*. For short, the sets $s_1, \ldots, s_k$ will be denoted $S^1$ and $s_{k+1} \ldots s_n$ with the symbol $S^2$; we will also write $s_i^1$ and $s_j^2$ to say that $s_i$ occurs in $S^1$ and $s_j$ in $S^2$. We allow $S$ to have either the first or the second component empty, i.e. we allow structured sets to have the form $\langle S^1 | \, \rangle$ or $\langle \, | S^2 \rangle$. Notice that we intend $\langle s_1 | s_2 \rangle$ to be different from $\langle s_2 | s_1 \rangle$, and $\langle s_1, s_2 | s_3 \rangle$ to be the same as $\langle s_2, s_1 | s_3 \rangle$.

A *space* $\mathcal{S}$ is a finite set of structured sets.

A structured set $S$ is called *trivial* when two sets of points $s_i, s_j$ are contained in different components of $S$ and are embeddable [3] one into the other (denoted $s_i \hookrightarrow s_j$) with respect to some prescribed class of maps (e.g. one can be interested to call trivial those structured sets where the embeddings are one-to-one maps between sets with the same structure. In case we consider tessellations $T$ for instance, we can think of the embeddings to be isometries (i.e. distance-preserving maps) of $T$ onto itself (faces onto faces and edges onto edges). Anything more complicated than this can be considered though.) We assume the *identity* map and the *empty embedding* (i.e. the empty set is always embeddable in any other set) always to be contained in any prescribed class of maps.

Notice that the components of a trivial structured set are non-empty. A space is called *trivial* if all its structured sets are trivial. The word *trivial* may seem very strong here, but in the motivating example of logic it is appropriate.

An *operator* of arity $l$ and subarity $n$ is a pseudomap from the space of $l$-tuples of structured sets into the space of structured sets such that

1. if $l$ structured sets $S_1 \ldots S_l$ are mapped by the operator to the structured set $S_0$, then there exist $n$ sets of points, say

$$
\begin{array}{cccccccc}
s_{1,1}^1, & \ldots & ,s_{j_1^1,1}^1 & \in & S_1^1 & \quad s_{1,1}^2, & \ldots ,s_{j_1^2,1}^2 & \in S_1^2 \\
\vdots & \vdots & \vdots & & & \vdots & \vdots & \vdots \\
s_{1,l}^1, & \ldots & ,s_{j_l^1,l}^1 & \in & S_l^1 & \quad s_{1,l}^2, & \ldots ,s_{j_l^2,l}^2 & \in S_l^2
\end{array}
$$

(where $j_1^1 + j_1^2 + \ldots + j_l^1 + j_l^2 = n$ and $j_i^1 + j_i^2 > 0$) and a set of points $s_0$ so that $S_0$ is defined to be either the structured set $\langle \{s_0\} \cup S_*^1 | S_*^2 \rangle$ or the structured set $\langle S_*^1 | \{s_0\} \cup S_*^2 \rangle$ where $S_*^1$ is

$$S_1^1 \setminus \{s_{1,1}^1, \ldots s_{j_1^1,1}^1\} \cup \ldots \cup S_l^1 \setminus \{s_{1,l}^1, \ldots s_{j_l^1,l}^1\}$$

and $S_*^2$ is

$$S_1^2 \setminus \{s_{1,1}^2, \ldots s_{j_1^2,1}^2\} \cup \ldots \cup S_l^2 \setminus \{s_{1,l}^2, \ldots s_{j_l^2,l}^2\}$$

---

[3] In this paper we call *embedding* any pseudomap from $X$ to $Y$ everywhere defined on $X$. An embedding is not required to preserve any structure on $X$ nor to be injective.

and where we ask each $s^i_{r,h}$ to be *embedded* into $s_0$, for $r = 1 \ldots j^i_h$, $h = 1 \ldots l$ and $i = 1, 2$ (once more there is no requirement on the embedding, beyond that it respect whatever structure has been imposed on the set $s_0$, we only ask for its existence.) All sets other than the $s_{r,h}$'s are embedded by the identity into the copy of themselves in $S_0$; and,

2. if $S'_1 \ldots S'_l$ are $l$ structured sets such that the same sets of points $s^i_{r,j}$ are contained in the $S'_j$'s in the same manner as above, then the operator maps $S'_1 \ldots S'_l$ into a structured set $S'_0$ in exactly the same manner as above.

The sets $s^i_{r,h}$ are called *arguments* of the operator, $s_0$ its *value* and $S_0$ its *output*.

The following figure illustrates a possible embedding of two finite graphs into a third one (we will only name those vertices needed to see how the embedding goes)

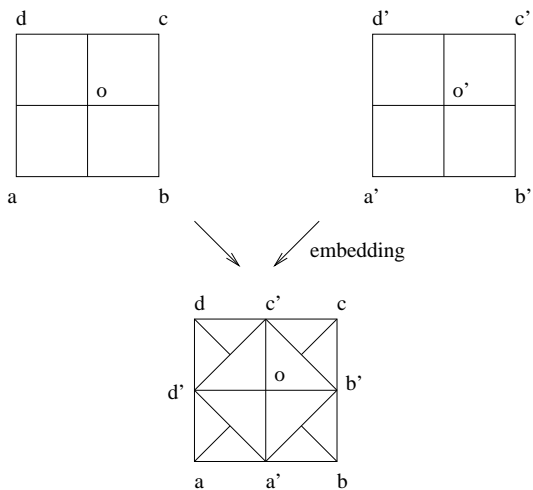

where the point $o'$ is mapped into $a, b, c, d$ by the embedding. This is allowed because an embedding is *not* defined as a map but as a pseudomap.

We should remark here that the number of arguments picked by the operator in the components of the structured sets on which it acts (i.e. $j^1_1, j^2_1, \ldots, j^1_l, j^2_l$) is fixed for each operator. Notice that we do not require operators to be defined on *all* families of structured sets and moreover, families of structured sets are not required to be ordered.

An operator will be called *surjective* when every point in its values $s_0$ lies in the image of at least one of the corresponding arguments $s^i_{r,h}$'s under the embeddings described above.

We say that a $l$-ary operator is *derived* from a $k$-ary operator (having larger subarity and $k \geq l$) if the following is true. Suppose that the $l$-ary operator acts on structured sets $S_1, \ldots, S_l$. Then we require that there be structured sets $S'_1, \ldots, S'_k$ on which the $k$-ary operator acts such that $S^1_j \subseteq {S'_j}^1$, $S^2_j \subseteq {S'_j}^2$ (for each $j = 1 \ldots l$), the $S'_j$'s contain only empty sets when $j = l+1 \ldots k$, and the arguments of the $k$-ary operator consist exactly of the arguments for the $l$-ary operator together with all the sets of points added to the $S_j$'s to get the $S'_j$'s when $j = 1\ldots l$, and all the empty sets in the $S'_j$'s when $j = l+1\ldots k$. Finally we require that the values and outputs of the two operators correspond in the obvious manner.

Whenever the auxiliary sets occurring in the $S'_j$'s (for $j = 1\ldots l$) are non-empty, the value of the derived operator is a set with a structure which 'depend' on the structure of the auxiliary sets (this follows directly from the definition of operator.) We will bring up this point again in Remark 3.5 and we will see that this dependence is conceptually innocuous.

A *transformation* on a space $\mathcal{S} = \{S_1, \ldots, S_k\}$ is induced by an operator acting on structured sets $S_{i_1}, \ldots, S_{i_l} \in \mathcal{S}$ and giving as a result the space $\mathcal{S}' = \{S_0\} \cup \mathcal{S} \setminus \{S_{i_1}, \ldots, S_{i_l}\}$ where $S_0$ is the output of the operator acting on $S_{i_1}, \ldots, S_{i_l}$; we ask all sets in $\mathcal{S}$ other than the ones involved in the action of the operator (i.e. the $s_{r,h}$'s in the definition of operator) to be embedded by the identity into the copy of themselves in $\mathcal{S}'$.

A pair of ($l$-ary and unary) operators acting on structured sets $S_1, \ldots, S_l$ and $S_r$ respectively, is called *regular* if

1. the $l$-ary operator takes exactly one argument in each $S_h$, say in the $i_h$-th component of $S_h$ (with $i_h \in \{1, 2\}$ and $h = 1, \ldots, l$), and

2. if $s_1 \ldots s_l$ are arguments for the $l$-ary operator acting on $S_1 \ldots S_l$, then $s_1 \ldots s_l$ also appear in $S_r$, but in the opposite components; these $s_h$'s are the set of arguments of the unary operator, and

3. the analogues of the set $s_0$ (in the above discussion of operators) for the $l$-ary operator and the unary operator are the same, and they belong to opposite components. The embeddings of the $s_h$'s into this common value $s_0$ are also the same.

Two operators belonging to a regular pair will be called *duals*. From the definition of regular pair of operators we can observe that the unary operator is acting in a space which is a sort of 'complement' to the space where the $l$-ary operator acts. For instance, suppose that our structured sets are graphs whose vertices are sets of points and they are colored either blue or red (we can always think of a bipartite collection equipped with some extra structure.) Let $S^1$ be the set of blue vertices and $S^2$ be the set of red ones, for any such graph $S$. Then, for any unary regular operator acting on blue vertices, the $l$-ary regular dual should act on red ones only.

We say that an operator is *constructible* from a set $\mathcal{O}$ of operators if either it belongs to $\mathcal{O}$ or it is obtained as a composition of operators in $\mathcal{O}$.

A set of operators $\mathcal{O}$ is *closed under regularity* if the following is true. Suppose that we are given an $l$-ary operator in $\mathcal{O}$ which acts on the structured sets $S_1, \ldots, S_l$, and let a sequence $i_0, i_1, \ldots, i_l \in \{1, 2\}$ be given. Suppose also that we have chosen sets of points $s_h \in S_h$ for each $h = 1\ldots l$, and that in fact each $s_h$ belongs to the $i_h$ component of $S_h$. Assume finally that we have another structured set $S_r$, not necessarily related to the $S_h$'s, which also contains each $s_h$, but in the opposite component from $i_h$. Then we ask that there exist a regular pair constructible from $\mathcal{O}$ such that the $l$-ary operator acts on the structured sets $S_1, \ldots, S_l$ taking a value in the $i_0$ component, the unary operator acts on $S_r$, and the conditions in the definition of regular pairs hold with the choice of $s_h$'s. (Note that the $l$-ary operator in the regular pair need not be the $l$-ary operator that we started with.)

This condition may seem complex but it is very natural in the motivating example of logic, and in any case it is just a way to combine some sets in a coherent manner.

We will call an operator *regular* if it belongs to a regular pair promised in the definition of closed under regularity. We will refer to it as $\mathcal{O}$-*regular* when we need to trace its dependence on a set of operators $\mathcal{O}$.

In the sequel, we will be interested in studying how trivial structured sets evolve through transformations induced by a set of operators $\mathcal{O}$ which is *closed* under regularity. We will say that a space $\mathcal{S}$ (equivalently its structured sets) is $\mathcal{O}$-*definable* if it is obtained by transforming a trivial space $\mathcal{S}'$ with (a finite number of applications of) operators which either are constructible from $\mathcal{O}$ or are derived from $\mathcal{O}$-regular operators. If $\mathcal{S}$ is $\mathcal{O}$-definable from $\mathcal{S}'$, it follows readily by definition that all points in $\mathcal{S}'$ are mapped into $\mathcal{S}$. In particular, because trivial structured sets require themselves the existence of an embedding, we will be thinking of points in $\mathcal{S}$ as being mapped into points of $\mathcal{S}$ in the obvious way (i.e. by composition of the embeddings or their inverses.) Not all points of $\mathcal{S}$ are expected to be mapped into $\mathcal{S}$ though. The composition of the embeddings and their inverses defines indeed a *pseudomap* between points in $\mathcal{S}$.

## 3 Interpolation for spaces of structured sets

Roughly, we will show that if $\mathcal{O}$ is a set of operators closed under regularity and defining a pseudomap $f$ between two sets $A$ and $B$ (from trivial spaces), then we can find a set $I$, a surjective pseudomap $f_A$ from $A$ to $I$ and an everywhere defined pseudomap $f_B$ from $I$ to $B$ so that for all $z \in I$ there is an $x \in A$ and a $y \in B$ such that $z \in f_A(x), y \in f_B(z)$ and $y \in f(x)$.

**3.1. Theorem.** (Interpolation Theorem for Structured Sets) *Let $S = \langle A|B \rangle$ be a structured set obtained from some trivial space $\mathcal{S}'$ by applying a finite number of transformations based on a set of operators $\mathcal{O}$ which is closed under regularity. Then there are two structured sets $S^A = \langle A|I \rangle$ and $S^B = \langle I|B \rangle$ (or $S^A = \langle A, I| \rangle$ and $S^B = \langle |B, I \rangle$) (where $I$ is called the* interpolant set *of $A, B$) which are both $\mathcal{O}$-definable. Moreover, if the regular operators are surjective operators and the trivial structured sets in $\mathcal{S}'$ are defined by surjective embeddings, then for all points $z \in I$ there are points $x \in A$ and $y \in B$ such that $z$ is mapped into $x$ and into $y$ by the transformations defining $S^A$ and $S^B$, and such that the points $x$ and $y$ are mapped into each other by the transformations defining $S$.*

*If $I$ does not contain any point (i.e. $I$ is built out of empty sets only) and no structured set in $\mathcal{S}'$ is trivial because of empty embeddings, then either $S^A = \langle A| \ \rangle$ or $S^B = \langle \ |B\rangle$ can be built from trivial spaces using only operators in $\mathcal{O}$.*

The ambiguity between $S^A = \langle A|I\rangle$, $S^B = \langle I|B\rangle$ versus $S^A = \langle A, I| \ \rangle$, $S^B = \langle \ |B, I\rangle$ may seem unsettling at first, but it is a natural symmetry in the problem. When the theorem will be interpreted in the context of logic for instance, this ambiguity will be naturally resolved due to the presence of rules for negation which allow one to move a formula from one side of a sequent to the other.

The interpolant set $I$ for the structured set $S$ will be built in steps from the definition of $S$ by applying either regular operators or operators derived from regular ones. If the regular operators are surjective operators and the trivial space is defined through surjective mappings, the construction will map a point in $I$ to a point in $A$ ($B$) first by using the inverses of the surjective embeddings to get back to the original trivial space, and then the embeddings used to define $A$ ($B$), to get into $A$ ($B$).

**Proof.** (Theorem 3.1) The proof is by induction on the steps of transformation from the trivial space $\mathcal{S}_0 = \{S_{0,1}, \ldots S_{0,h}\}$ to the space $\mathcal{S} = \{S\}$. Call $\mathcal{S}_{p+1}$ the space obtained from $\mathcal{S}_p$ by the $p+1$-th step of transformation from $\mathcal{S}_0$ to $\mathcal{S}$.

Before we begin the proof in earnest let us consider a simplified argument which is not strong enough for Theorem 3.1 but which brings out some of the key points. Notice that if $s$ is a set of points which belongs to a structured set in $\mathcal{S}_p$ (for some $p$), then the transformations which convert $\mathcal{S}_p$ into $\mathcal{S}$ will induce an embedding of $s$ into either $A$ or $B$. Thus for each $p$ we can consider the spaces $\mathcal{S}_p(A)$ and $\mathcal{S}_p(B)$ which are obtained by taking the structured sets in $\mathcal{S}_p$ and removing from them the sets of points which are embedded into $B$ or $A$, respectively. This simple-minded decomposition has the nice feature that the transformation that takes $\mathcal{S}_p$ to $\mathcal{S}_{p+1}$ induces similar transformations for $\mathcal{S}_p(A)$ and $\mathcal{S}_p(B)$. This is not good enough for the theorem though, because $\mathcal{S}_0(A)$ and $\mathcal{S}_0(B)$ are typically not trivial spaces. Although the transformations do not mingle the sets embedded into $A$ with the sets embedded into $B$, the requirement that the initial spaces be trivial necessitates some links which are more subtle. The idea of this trivial model is clearly present in the induction argument below, but it is necessary to also choose carefully some interpolating sets, starting at level $p = 0$.

Our induction argument goes as follows. At step $p$ we consider the space $\mathcal{S}_p = \{S_{p,1} \ldots S_{p,k}\}$ and we show that

1. there are two spaces $\mathcal{S}_p^A = \{S_{p,1}^A \ldots S_{p,i_r}^A\}$ and $\mathcal{S}_p^B = \{S_{p,1}^B \ldots S_{p,i_l}^B\}$ (with $r, l \leq k$) such that

    (a) $\mathcal{S}_p^A$ and $\mathcal{S}_p^B$ are trivial spaces when $p = 0$ and they are $\mathcal{O}$-definable [4] (from $\mathcal{S}_0^A$ and $\mathcal{S}_0^B$) when $p > 0$, and

---

[4] Notice that we are permitting the use of operators derived from regular operators in $\mathcal{O}$-definability.

    (b) each structured set in $\mathcal{S}_p^A$ is associated to a structured set in $\mathcal{S}_p$, and for each structured set in $\mathcal{S}_p$ there is at most one counterpart in $\mathcal{S}_p^A$; the same is true for $\mathcal{S}_p^B$; each structured set in $\mathcal{S}_p$ has a counterpart in at least one of $\mathcal{S}_p^A$, $\mathcal{S}_p^B$, and

    (c) if $S_{p,i}$ is a structured set in $\mathcal{S}_p$, and if $S_{p,i}$ contains sets which are embedded into $A$, then $S_{p,i}$ is associated to a structured set in $\mathcal{S}_p^A$ which contains these same sets, and in the same components; if $S_{p,i}$ is a structured set in $\mathcal{S}_p$ which is associated to a structured set in $\mathcal{S}_p^A$, then the latter contains *at most* one other set of points beyond the ones just mentioned (an $I$-set, see below); the same is true for $\mathcal{S}_p^B$, and

2. there exist special sets called $I$-sets which will be used to build the interpolant set $I$ such that

    (a) if $S_{p,i}$ is a structured set in $\mathcal{S}_p$, and if it is associated to structured sets in each of $\mathcal{S}_p^A$, $\mathcal{S}_p^B$, then each of these structured sets has an $I$-set, the two $I$-sets are the same sets, and these $I$-sets lie in opposite components; if there is no structured set associated to $S_{p,i}$ in one of $\mathcal{S}_p^A$, $\mathcal{S}_p^B$ then the structured set in the other does not have an $I$-set, and

    (b) if the regular operators are surjective operators and the trivial structured sets in $\mathcal{S}_0$ are defined by surjective embeddings, then for any two structured sets $S_{p,i_j}^A$ and $S_{p,i_s}^B$ (for some $j, s$) associated to the same $S_{p,i}$ in $\mathcal{S}_p$, and all points $z$ in their $I$-set there are points $x$ in (some set of) $S_{p,i_j}^A$ and $y$ in (some set of) $S_{p,i_s}^B$ (for some $j, s$) such that $z$ is mapped into $x$ and into $y$ by those embeddings and their inverses defining $\mathcal{S}_p^A$ and $\mathcal{S}_p^B$, and such that the points $x$ and $y$ are mapped into each other by the transformations defining $S$.

We begin by defining the trivial spaces $\mathcal{S}_0^A = \{S_{0,1}^A \ldots S_{0,m}^A\}, \mathcal{S}_0^B = \{S_{0,1}^B \ldots S_{0,n}^B\}$ (with $m, n \leq h$) from which $\mathcal{S}^A = \{S^A\}$ and $\mathcal{S}^B = \{S^B\}$ will be obtained.

We will consider all the structured sets $S_{0,i}$ in $\mathcal{S}_0$, which are all trivial since $\mathcal{S}_0$ is assumed to be a trivial space, and we will define (whenever possible) $S_{0,i}^A$ in $\mathcal{S}_0^A$ and $S_{0,i}^B$ in $\mathcal{S}_0^B$ out of $S_{0,i}$'s together with a special set $I_{0,i}$.

Let $s_1, s_2$ be the embeddable sets which belong to the trivial structured set $S_{0,i}$. Without loss of generality suppose $s_1 \in S_{0,i}^1$ and $s_2 \in S_{0,i}^2$. There are four cases to consider:

1. All points in $s_1, s_2$ are embedded in $A$ of $\mathcal{S}$ (through the transformations that take $\mathcal{S}_0$ to $\mathcal{S}$). Call $s_{0,i_1}^A \ldots s_{0,i_l}^A$ those sets in $S_{0,i}$ embedded into $A$ of $\mathcal{S}$, and $s_{0,i_1}^B \ldots s_{0,i_r}^B$ those embedded into $B$ of $\mathcal{S}$, where $l + r$ is the number of sets in $S_{0,i}$. By using the division into components inherited from $S_{0,i}$, let

$$S_{0,i}^A = \langle s_{0,i_1}^A \ldots s_{0,i_q}^A | s_{0,i_{q+1}}^A \ldots s_{0,i_l}^A, \{\}\rangle$$

so that $S^A_{0,i}$ is trivial because of the embedding between $s_1, s_2$ that comes from the triviality of $S_{0,i}$ (notice that $s_1, s_2$ belong to distinct partition sets of $S^A_{0,i}$), and let

$$S^B_{0,i} = \langle \{\}, s^B_{0,i_1} \ldots s^B_{0,i_v} | s^B_{0,i_{v+1}} \ldots s^B_{0,i_r} \rangle$$

with a trivial embedding between the empty set $\{\}$ and say the set $s^B_{0,i_{v+1}}$, if $r > v$. Let $I_{0,i}$ in $S^A_{0,i}$ and $S^B_{0,i}$ be the set $\{\}$.

If $r = v$ and $v = 0$, let $S^A_{0,i} = S_{0,i}$ and $S^B_{0,i}$ be not defined. No $I$-set will be defined when $r = v = 0$.

If $r = v$ and $v > 0$, let $S^A_{0,i} = \langle \{\}, s^A_{0,i_1} \ldots s^A_{0,i_q} | s^A_{0,i_{q+1}} \ldots s^A_{0,i_l} \rangle$ with the embedding between $s_1, s_2$ as above, and $S^B_{0,i} = \langle s^B_{0,i_1} \ldots s^B_{0,i_v} | \{\} \rangle$ with a trivial embedding between $\{\}$ and say the set $s^B_{0,i_1}$. Let $I_{0,i}$ in $S^A_{0,i}$ and $S^B_{0,i}$ be the set $\{\}$.

2. All points in $s_1, s_2$ are embedded in $B$ of $\mathcal{S}$. This case is similar to 1.

3. All points in $s_1$ are embedded in $A$ of $\mathcal{S}$ and all points in $s_2$ are embedded in $B$ of $\mathcal{S}$. Suppose $s_1 \hookrightarrow s_2$ ($s_2 \hookrightarrow s_1$.) By using the division into components inherited from $S_{0,i}$, define

$$S^A_{0,i} = \langle s^A_{0,i_1} \ldots s^A_{0,i_q} | s^A_{0,i_{q+1}} \ldots s^A_{0,i_l}, s_1 \rangle$$
$$(S^B_{0,i} = \langle s_2, s^B_{0,i_1} \ldots s^B_{0,i_v} | s^B_{0,i_{v+1}} \ldots s^B_{0,i_r} \rangle)$$

with the identity map as an embedding of $s_1$ ($s_2$) into $s_1$ ($s_2$) (notice that $s^A_{0,i_j}$ ($s^B_{0,i_j}$) is $s_1$ ($s_2$) for some $j = 1 \ldots q$ ($j = v+1 \ldots r$)), and

$$S^B_{0,i} = \langle s_1, s^B_{0,i_1} \ldots s^B_{0,i_v} | s^B_{0,i_{v+1}} \ldots s^B_{0,i_r} \rangle$$
$$(S^A_{0,i} = \langle s^A_{0,i_1} \ldots s^A_{0,i_q} | s^A_{0,i_{q+1}} \ldots s^A_{0,i_l}, s_2 \rangle)$$

with the triviality coming from the embedding between $s_1, s_2$ used to define $S_{0,i}$ as a trivial space (notice that $s^B_{0,i_j}$ ($s^A_{0,i_j}$) is $s_2$ ($s_1$) for some $j = v+1 \ldots r$ ($j = 1 \ldots q$).) Let the space $I_{0,i}$ in $S^A_{0,i}$ and $S^B_{0,i}$ be the set $s_1$ ($s_2$.)

4. All points in $s_1$ are embedded in $B$ of $\mathcal{S}$ and all points in $s_2$ are embedded in $A$ of $\mathcal{S}$. This case is symmetric to 3.

By making these constructions for each $S_{0,i}$ we obtain the trivial spaces $\mathcal{S}^A_0 = \{S^A_{0,1}, \ldots, S^A_{0,m}\}$, $\mathcal{S}^B_0 = \{S^B_{0,1}, \ldots, S^B_{0,n}\}$ as desired, and we also obtain (whenever defined) $I$-sets $I_{0,i}$ in $S^A_{0,i}$ and $S^B_{0,i}$. Next we want to define spaces $\mathcal{S}^A_p$ and $\mathcal{S}^B_p$ for each $p$, with $p = 0$ corresponding to the construction just made, and we also want to define $I$-sets for these spaces.

Suppose that we have made our construction to stage $p$ and we want to go to stage $p + 1$. Suppose also that $\mathcal{S}_{p+1}$ is obtained from $\mathcal{S}_p$ by applying a $k$-ary operator to structured sets $S_1, \ldots, S_k$ in $\mathcal{S}_p$. Let $s_0$ be the set of points in the output of this $k$-ary operator as described in the definition of transformations in section 2, so that $s_0$ is an element of one of the structured spaces in $\mathcal{S}_{p+1}$. Assume that $s_0$ embeds into $A$ by the mappings that are implicit in the transformations which convert $\mathcal{S}_{p+1}$ into $\mathcal{S} = \{S\} = \{\langle A|B \rangle\}$. (Otherwise it embeds into $B$ and one employs a similar construction.) We want to define $\mathcal{S}_{p+1}^A$ and $\mathcal{S}_{p+1}^B$. Our induction hypothesis provides us with spaces $\mathcal{S}_p^A$ and $\mathcal{S}_p^B$ and structured spaces $S_1^A, \ldots, S_k^A$ in $\mathcal{S}_p^A$, $S_1^B, \ldots, S_l^B$ in $\mathcal{S}_p^B$ (for $0 \leq l \leq k$; the inequality is a consequence of the possible undefinability of structured sets in cases 1 and 2 for $p = 0$) which correspond to $S_1, \ldots, S_k$ in $\mathcal{S}_p$. To define $\mathcal{S}_{p+1}^A$ we basically want to apply the same $k$-ary operator to $S_1^A, \ldots, S_k^A$ as was applied to $S_1, \ldots, S_k$, and roughly to replace the structured spaces $S_1^B, \ldots, S_l^B$ by their union to get $\mathcal{S}_{p+1}^B$, but this is not quite right, because we would get too many $I$-sets in the structured space that results. To fix this we use the assumption of $\mathcal{O}$ being closed under regularity.

Formally, apply the $k$-ary operator to the sets $S_1^A, \ldots S_k^A$ and call $\mathcal{S}_{p+1,*}^A$ the space that results from $\mathcal{S}_p^A$. Notice that the resulting structured set $S_0$ in $\mathcal{S}_{p+1,*}^A$ will contain $l \geq 0$ $I$-sets coming from $S_1^A, \ldots S_k^A$. If $l = 0, 1$ let $\mathcal{S}_{p+1}^A$ be $\mathcal{S}_{p+1,*}^A$. If $l > 1$ then apply a suitable regular *unary* operator to the $l$ $I$-sets (in case $l = k$, apply a unary regular operator of subarity $k$; in case there are $l < k$ $I$-sets among $S_1^A, \ldots S_k^A$, consider a unary operator of subarity $l$ derived from the unary operator of subarity $k$ by adding $k - l$ empty sets to the structured set $S_0$ in $\mathcal{S}_{p+1,*}^A$, as described in the definition of derived operator in section 2) and call $\mathcal{S}_{p+1}^A$ the space that results from $\mathcal{S}_{p+1,*}^A$. Notice that such a regular unary operator of subarity $k$ exists because our system is closed under regularity, i.e. for any $k$-ary operator in $\mathcal{O}$ and for any $k_1, k_2$ such that $k_1 + k_2 = k$ there is a regular unary operator that when applied to $\langle S^1|S^2 \rangle$ takes $k_1$ arguments in $S^1$ and $k_2$ arguments in $S^2$.

We would like to apply now the dual $k$-ary regular operator to the $I$-sets in $S_1^B, \ldots S_l^B$ of $\mathcal{S}_p^B$. First we rule out the cases $l = 0, 1$ by taking $\mathcal{S}_{p+1}^B$ to be $\mathcal{S}_p^B$; if $l = k$ (for $k > 1$) then we apply the dual $k$-*ary* regular operator as expected, otherwise (i.e. $k - l > 0$ and $l > 1$) we apply a $l$-ary operator which we derive from the $k$-ary dual regular operator by considering $k - l$ of its arguments to be structured sets containing only the empty set occurring in the appropriate component (the component is a priori determined because the operator is dual to the one used to define $\mathcal{S}_{p+1}^A$.) The resulting set is the $I$-set of the space $\mathcal{S}_{p+1}^B$; it is the same set as the one produced for $\mathcal{S}_{p+1}^A$.

In its final stage, we will have $\mathcal{S}_p = \mathcal{S} = \{\langle A \mid B \rangle\}$, and $\mathcal{S}_p^A, \mathcal{S}_p^B$ spaces which each consists of only one structured set. The structured set for $\mathcal{S}_p^A$ consists of $A$ and an interpolant set $I$, and the structured set for $\mathcal{S}_p^B$ consists of $B$ and the same interpolant set. Since $I$-sets are built only out of dual operators and dual operators take their values in opposite components, at the end of the construction we have either built a pair of structured sets $\langle A \mid I \rangle$, $\langle I \mid B \rangle$ or two structured sets $\langle A, I \mid \rangle$, $\langle \mid B, I \rangle$.

The assertion of the theorem concerning surjective regular operators and trivial structured sets in $\mathcal{S}_0$ defined by surjective embeddings, corresponds to condition 2.(b) in the induction hypothesis, for $p$ being the final stage. Notice that the condition is satisfied at all steps $p$ since the derived operators we employ themselves use only auxiliary sets which are empty. This together with the assumption that the regular operators are surjective implies the property we need.

Let us now prove the last assertion of the theorem, that if $I$ does not contain any points, and if the triviality of the structured sets in the original trivial space $\mathcal{S}_0$ never relies on empty embeddings, then either $\langle A|\,\rangle$ or $\langle\,|B\rangle$ is already constructible from a trivial space using only operators in $\mathcal{O}$.

By assumption $\mathcal{S} = \{\langle A|B\rangle\}$ can be derived from the trivial space $\mathcal{S}_0$ by finitely many transformations in $\mathcal{O}$, and as before we let $\mathcal{S}_p$ denote the spaces obtained in the intermediate steps.

Given a set of points $s$ in a structured set in some $\mathcal{S}_p$, we call $s$ an *A-set* (*B-set*) if it is eventually mapped into $A$ ($B$) by the transformations which convert $\mathcal{S}_0$ into $\mathcal{S}$. Of course one of the two possibilities always obtains.

Our assumption that $I$ contains no points implies that each structured set $S$ in $\mathcal{S}_0$ is trivial either because of a pair of $A$-sets $s_1, s_2$ in opposite components of $S$ with $s_1$ embeddable into $s_2$ (in which case we say that $S$ is *A-trivial*), or because of a pair of $B$-sets with the same properties ($S$ is *B-trivial*). That is, the triviality of $S$ cannot require an embedding between an $A$-set and a $B$-set, because an embedding between an $A$-set and a $B$-set would have to be trivial (empty) since the interpolant $I$ has no elements. This follows from our earlier construction. Note that this is the only way in which we use the hypothesis that $I$ has no elements, and we do not need to assume that $\mathcal{O}$ is closed under regularity for the argument that follows.

In order to show that $\langle A|\,\rangle$ or $\langle\,|B\rangle$ can be obtained from a trivial space through transformations only induced by operators in $\mathcal{O}$ we want to use the assumption that $\mathcal{S}$ is a result of a finite number of transformations based on $\mathcal{O}$ and we want to build the new derivation by simplifying the old one. The main point of the simplification lies already in the observation that if an $A$-trivial set is acted on by a $B$-transformation, then the result is $A$-trivial and so the transformation is not really needed. To proceed we may as well assume that our given derivation of $\mathcal{S}$ from $\mathcal{S}_0$ by transformations in $\mathcal{O}$ cannot be further simplified, i.e. it has the smallest number of steps possible in the sense that there is no subsequence $\mathcal{S}_0 \ldots \mathcal{S}_{j+1}$ of $\mathcal{S}_0 \ldots \mathcal{S}_p = \mathcal{S}$ with $j < p$ that can be redefined as a sequence $\mathcal{S}'_0 \ldots \mathcal{S}'_j$ such that $\mathcal{S}'_j = \mathcal{S}_{j+1}$ and such that for each $i = 0 \ldots j-1$ the space $\mathcal{S}'_{i+1}$ is obtained from $\mathcal{S}'_i$ using the same operator applied to $\mathcal{S}_i$ to obtain $\mathcal{S}_{i+1}$ (with the same sets of points as arguments and values.)

Let $S$ be a structured set in some $\mathcal{S}_p$. The structured sets in $\mathcal{S}_0$ which are sent to $S$ by the transformations that convert $\mathcal{S}_0$ into $\mathcal{S}_p$ are called the *ancestors* of $S$ in $\mathcal{S}_0$. We say that $S$ is *A-justified* if every ancestor of $S$ in $\mathcal{S}_0$ is $A$-trivial, and we define *B-justified* in the same way.

We claim that every structured set $S$ in any $\mathcal{S}_p$ is either $A$-justified or $B$-justified. To see this we argue by induction. Our assumptions guarantee that this is true when $p = 0$. Suppose that it is true at level $p$ and that we want to prove it at level $p+1$. Let $S$ be a structured set in $\mathcal{S}_{p+1}$. Then either $S$ was already present in $\mathcal{S}_p$, in which case there is nothing to prove, or $S$ arises as the output of some structured sets $S_1, \ldots, S_k$ in $\mathcal{S}_p$ under a $k$-ary transformation in $\mathcal{O}$. Let $s_0$ be the set of points which is the value of this transformation, as in section 2. Assume that $s_0$ is an $A$-set, the argument is the same if it is a $B$-set. Each of the $S_i$'s is either $A$-justified or $B$-justified, by induction hypothesis, and they cannot all be one or the other, since $S$ itself would then be too. So the $S_i$'s cannot all be $A$-justified, at least one is $B$-justified, let us say $S_1$.

Consider the $B$-part of $S_1$, i.e. the structured set which you get by keeping the $B$-sets in $S_1$ and throwing away the $A$-sets. Since $s_0$ is assumed to be an $A$-set we have that $S$ in $\mathcal{S}_{p+1}$ must contain an exact copy of the $B$-part of $S_1$ inside of itself, in fact nothing in there can be mapped into $s_0$.

Notice now that $S_1$ is obtained from its ancestors in $\mathcal{S}_0$ in $n \leq p$ transformations. (It takes $p$ steps to convert $\mathcal{S}_0$ into $\mathcal{S}_p$, but some of these steps may not affect the ancestors of $S_1$, and we only want to count those steps.) The $B$-part of $S_1$ can be obtained from the $B$-part of the ancestors of $S_1$ in $\mathcal{S}_0$ by $\leq n$ of these transformations. Indeed, we simply use the transformations for obtaining $S_1$ that operate on the $B$-parts and we throw away the transformations that affect the $A$-parts. Since the ancestors of $S_1$ in $\mathcal{S}_0$ are $B$-trivial by assumption we conclude that the $B$-part of $S_1$ is obtained from a trivial space with $\leq n$ steps of transformations based on $\mathcal{O}$.

This implies that $S$ can be obtained in $\leq n$ steps from trivial spaces, because $S$ contains the $B$-part of $S_1$ and we can always add as many sets of points as we want (i.e. we can add them to the trivial spaces and then drag them along).

This contradicts our assumption of minimality, because it provides a way for us to obtain $S$ in $\leq n \leq p$ steps, whereas it should take $p+1$. (Roughly speaking, we eliminated the need for the $p$-th step of transformation. Note that in order to contradict minimality we make changes to the previous steps only when they involve the ancestors of $S$.)

This completes the proof of our claim that every structured set in every $\mathcal{S}_p$ is either $A$-justified or $B$-justified.

Now we apply the claim to the structured set $\langle A|B\rangle$ in the final space. We conclude that it is either $A$-justified or $B$-justified, let us say $A$-justified (the other case is the same). For the same reason as before we conclude that the $A$-part of this structured set, namely $\langle A|\,\rangle$, is obtained from a trivial space through transformations based on $\mathcal{O}$. This is exactly what we wanted.

Note that this use of minimality may seem nonconstructive, but it is easy to convert the preceding argument into one which is constructive (but then less pleasant to read). Also we could have derived more information from the assumption of minimality, but it was not needed and so we did not bother.

This completes the proof of Theorem 3.1. ∎

**3.2. Remark.** Our definition of a set of operators being *closed under regularity* requires regular pairs to satisfy a *minimal* set of conditions sufficient to prove Theorem 3.1. Because of this generality, the number of regular pairs required to belong to such a set of operators might be quite large. Theorem 3.1 can be carried out for a restricted number of regular pairs whenever constraints on the operators defining a space are given. For instance if we can define a space $\mathcal{S} = \{\langle A|B\rangle\}$ through operators whose arguments and values both belong to the same component (a fixed one for each operator) we can ask to the $I$-sets in the space $\mathcal{S}_p^A$ ($\mathcal{S}_p^B$) of our construction, to always belong to the second (first) component of its structured sets. Therefore we would only need to consider regular $l$-ary operators whose arguments and values lie both either in the first or in the second component of structured sets.

**3.3. Remark.** We say that $S$ is a *restricted* structured set if its second component $s_{k+1}\ldots s_n$ contains only copies of the same set (where each copy should be equipped with the same structure), i.e. $s_{k+1} = \ldots = s_n$. Let us consider spaces and operators defined with respect to restricted structured sets. A *restricted transformation* on a space $\mathcal{S}$ is a transformation induced by an operator acting only on restricted structured sets and giving as output a restricted structured set. We say that a space $\mathcal{S}$ is obtained by *restricted $\mathcal{O}$-definability* if it is $\mathcal{O}$-definable through restricted transformations. A version of Theorem 3.1 based on restricted $\mathcal{O}$-definability, can be proved for any set of operators $\mathcal{O}$ which is not only closed under regularity but also contains constructible $n$-ary operators mapping $n$ copies of the same set which lie on the $i$-th component into the set itself lying in the $i$-th component (for all $n \geq 1$ and $i \in \{1,2\}$; these latter operators will be called operators of *contraction* as their function suggests.) We say that a set of operators so defined is closed under *extended regularity* condition.

We will see in section 4 that intuitionistic logic is interpretable as a space of restricted structured sets. Theorem 3.4 will turn out to be the set-theoretical interpretation of the Craig Interpolation Theorem for the intuitionistic logical system $LJ$.

**3.4. Theorem.** (Interpolation Theorem for Restricted Structured Sets) *Let $\mathcal{O}$ be a set of operators which is closed under extended regularity. Let $\mathcal{S} = \{S\} = \{\langle A|B\rangle\}$ be a space obtained from some trivial space $\mathcal{S}'$ by applying a finite number of restricted transformations based on $\mathcal{O}$. Suppose that for each restricted structured set involved in the transformation of $\mathcal{S}'$ into $\mathcal{S}$, any pair of sets occurring in its second component is embedded either into $A$ or into $B$. Then there are two restricted structured sets $S^A = \langle A|I\rangle$ and $S^B = \langle I|B\rangle$ (for some set $I$ called interpolant of $A, B$) which are both obtained by restricted $\mathcal{O}$-definability. Moreover, if regular operators are surjective and the trivial structured sets in $\mathcal{S}'$ are defined by surjective embeddings, then for all points $z \in I$ there are points $x \in A$ and $y \in B$ such that $z$ is mapped into $x$ and into $y$ by the transformations defining $S^A, S^B$,* and such that the points $x$ and $y$ are mapped into each other by the transformations defining $\mathcal{S}$.

If $I$ does not contain any point (i.e. $I$ is built out of empty sets only) and no structured set in $\mathcal{S}'$ is trivial because of empty embeddings, then either $S^A = \langle A| \ \rangle$ or $S^B = \langle \ |B\rangle$ can be built from trivial spaces using only operators in $\mathcal{O}$.

**Proof.** Similar to the proof of Theorem 3.1. Before we begin to spell out the differences with the proof of Theorem 3.1 we want to remind the reader that any pair of sets occurring in the second component of any of the structured sets that we consider (which are involved in the transformation converting $\mathcal{S}'$ into $\mathcal{S}$), will be embedded (by the transformation) either into $A$ or into $B$. In the sequel we will frequently use this fact and we will not be explicit about it anymore.

As in the proof of Theorem 3.1 we will build $I$-sets for structured sets $S_i^A, S_i^B$ in such a way that they will lie in opposite components and they will not interfere with the restriction on the structured sets. In particular, at each stage $p$ we will show that if a structured set in $\mathcal{S}_p^A$ or $\mathcal{S}_p^B$ has an $I$-set in its second component, then that is the only set of points in the second component of that structured set. We will continue to use the same notation as in the proof of Theorem 3.1.

At stage $p = 0$, case 1 implies $r = v$. If $v = 0$ there will be no $I$-set defined; if $v > 0$ then the $I$-set is defined to be the empty set and it will be lying on the second component of $S_{0,i}^B$. Notice that the second component of $S_{0,i}^B$ will not contain other sets than the empty set since $s_2$ in $S_{0,i}$ is mapped into $A$. The same holds for case 2. In case 3, the $I$-set is defined to be $s_1$ ($s_2$) and we only need to observe that the structured sets $S_{0,i}^A = \langle s_{0,i_1}^A \ldots s_{0,i_q}^A | s_1 \rangle$ ($S_{0,i}^B = \langle s_2, s_{0,i_1}^B \ldots s_{0,i_v}^B | s_2, \ldots, s_2 \rangle$) and $S_{0,i}^B = \langle s_1, s_{0,i_1}^B \ldots s_{0,i_v}^B | s_2, \ldots, s_2 \rangle$ ($S_{0,i}^A = \langle s_{0,i_1}^A \ldots s_{0,i_q}^A | s_2 \rangle$) are restricted. The same holds for case 4.

At stage $p + 1$, in defining $\mathcal{S}_{p+1}^A, \mathcal{S}_{p+1}^B$ from $\mathcal{S}_p^A, \mathcal{S}_p^B$ we need to change the order of construction of the $I$-sets proposed in Theorem 3.1. Namely, we build first the $I$-sets and then the set $s_0$ (i.e. the value of the $k$-ary operator acting on $\mathcal{S}_p$ at stage $p + 1$) which we assume to be mapped into $A$ by the implicit transformation from $\mathcal{S}_{p+1}$ to $\mathcal{S}$ (the construction goes similarly if $s_0$ is embedded into $B$.)

Let us consider some cases.

First, suppose that $\mathcal{S}_{p+1}$ is obtained from $\mathcal{S}_p$ by applying a $k$-ary operator to $S_1 \ldots S_k$ in $\mathcal{S}_p$ which takes a value in the *second* component of its output $S_0$ in $\mathcal{S}_{p+1}$. Because of the restriction on structured sets that should be satisfied, we will not immediately apply the $k$-ary operator to $S_1^A, \ldots, S_k^A$ in $\mathcal{S}_p^A$; instead, for all $S_i^A$ containing an $I$-set we will first apply a unary operator of subarity 1 to the $I$-set in $S_i^A$ which should be *derived* from a regular unary operator of subarity $k$ by using $k - 1$ auxiliary sets, each of them being a copy of the $I$-sets in $S_1^A, \ldots, S_{i-1}^A, S_{i+1}^A, \ldots, S_k^A$, and adding them to $S_i^A$ (keeping the same component) as in the definition of derived unary operator in section 2 [5] (note that if there are $l < k$ $I$-sets in $S_1^A, \ldots, S_k^A$ then we should fix $k - l$ sets to be the empty set; we ask the empty sets to lie in the first component, but this choice is arbitrary); we will ask the derived regular unary operator to take its value in the *first* component of its output set. If $S_i^A$ does *not* contain an $I$-set we do not do anything (this can happen if there are $l < k$ $I$-sets.)

---

[5] The unary operator of subarity $k$ we consider here is not required to act on restricted structured sets. This is consistent with the definition of restricted $\mathcal{O}$-definability where we ask only that the derived operators act on restricted structured sets.

This construction should be repeated for all $i = 1 \ldots k$. Since all unary operators that we applied have been derived from the same regular unary operator, their values are copies of the same set $s_0$ and all lie in the first component (even though in different structured sets.) As a result of this operation we obtain a space with structured sets $S_1^{A,*} \ldots S_k^{A,*}$ defined as $S_1^A \ldots S_k^A$ except for the $I$-sets that are replaced by $s_0$.

We apply now the $k$-ary operator to $S_1^{A,*} \ldots S_k^{A,*}$ and obtain a structured set containing $l$ copies of $s_0$ lying in the first component. We want to reduce them to one set and for this we apply the operator of contraction in $\mathcal{O}$ which maps all copies of $s_0$ into one (this operator exists by the assumption of extended regularity.) Call $\mathcal{S}_{p+1}^A$ the space so obtained.

Apply the dual regular $k$-ary operator to the $I$-sets in $S_1^B \ldots S_l^B$ of $\mathcal{S}_p^B$ as indicated in Theorem 3.1 and call $\mathcal{S}_{p+1}^B$ the resulting space. Since the unary regular operators of subarity 1 used to define $\mathcal{S}_{p+1}^A$ were all derived by the same unary regular operator of subarity $k$ which was required to take a value in the first component, the *dual* $k$-ary operator will take a value in the second component. This choice will not be in conflict with the restriction on structured sets because of condition 1.(c) in the induction argument (see the proof of Theorem 3.1). In fact by hypothesis all sets in the second component of the output $S_0$ in $\mathcal{S}_{p+1}$ are embedded into $A$. This concludes the first case.

Now, suppose that a $k$-ary operator acts on structured sets $S_1 \ldots S_k$ in $\mathcal{S}_p$ and takes a value in the *first* component of a structured set $S_0$ of $\mathcal{S}_{p+1}$. The construction is similar to the case just discussed except for minor differences about where to put the $I$-set. If the second component of the structured set $S_0$ in $\mathcal{S}_{p+1}$ is not empty, and the sets contained in it are embedded into $A$ by the transformation converting $\mathcal{S}_{p+1}$ into $\mathcal{S}$, then the $I$-set in $\mathcal{S}_{p+1}^A$ should be required to lie in the first component. Otherwise (the sets are embedded into $B$ and) the $I$-sets in $\mathcal{S}_{p+1}^A$ should be required to lie in the second component. This can always be done because by hypothesis the set of operators $\mathcal{O}$ is closed under extended regularity, and therefore for any choice of component $i \in \{1, 2\}$ there are regular operators and operators of contraction whose values lie in the $i$-th component. The rest of the construction proceeds as in the first case.

We should point out here that if an $I$-set in $\mathcal{S}_{p+1}^A$ (for some $p$) is obtained from $\mathcal{S}_p^A$ through derived operators based on auxiliary sets which are copies of non empty $I$-sets, then for each of these $I$-sets there is always a structured set $S_i^A$ in $\mathcal{S}_p^A$ (for $i = 1 \ldots k$) which contains it and satisfies condition 2.(b) at stage $p$. This is enough to imply condition 2.(b) at stage $p+1$ since one can show that the transformation of $\mathcal{S}_p$ into $\mathcal{S}_{p+1}$ is surjective. ∎

**3.5. Remark.** While in the proof of Theorem 3.1 we need operators derived by adding only copies of the empty set, it is worthwhile to notice that in the proof of Theorem 3.4 operators are derived by using auxiliary sets which are in general non-trivial, i.e. they are $I$-sets equipped with a certain structure. This induces transformations based on derived operators which have as a value a set with a structure somehow 'dependent' on the structure of these auxiliary sets. Clearly empty sets do not impose any structural requirement. In this sense the reader can be lead to the impression that in the proof of Theorem 3.4 we use a substantially different construction than in the proof of Theorem 3.1. Indeed this is not the case since one can show that any structured set so obtained can be $\mathcal{O}$-defined from non-restricted structured sets. The basic point (and this is a general fact) is that for any $\mathcal{O}$-definable structured set $\langle S^1|S^2\rangle$ (obtained with the application of $p$ steps of transformation) and any $S^1_* \supseteq S^1, S^2_* \supseteq S^2$, one can $\mathcal{O}$-define the structured set $\langle S^1_*|S^2_*\rangle$ (with $p$ steps of transformation) by adding to the trivial structured sets defining $\langle S^1|S^2\rangle$ those sets of points in $S^1_*\backslash S^1, S^2_*\backslash S^2$ and then drag them along the transformations as required. Once we know this, we are done. In fact, suppose we can derive a structured set $\langle S^1_0|S^2_0\rangle$ from an $\mathcal{O}$-definable structured set $\langle S^1|S^2\rangle$ by applying a unary operator of subarity $l < k$ derived from a unary operator of subarity $k$ by using $k - l$ auxiliary sets. Then, we can notice that by 'adding' to the suitable components of $\langle S^1|S^2\rangle$ the desired $k - l$ sets (using the general fact pointed out before), one can obtain $\langle S^1_0|S^2_0\rangle$ by the direct application of the unary operator of subarity $k$. (What we just described has a well-known counterpart in the context of logic with the distinction between multiplicative and additive unary rules; notice that the basic property we used corresponds to the *weakening* rule (see section 4).) The same idea can be applied to derived operators of arity $k$ and subarity $l < k$ to show that $\langle S^1_0|S^2_0\rangle$ can be obtained by applying the $l$-ary derived operator only to auxiliary sets which are empty.

**3.6. Remark.** If we define a structured set simply as a *collection* of sets $s_1 \ldots s_n$ (equipped with some additional structure) instead of a bipartite collection, we can show an interpolation theorem for sets $S = \{A, B\}$, where regular pairs are the same as before except that we forget the bipartite structure.

**3.7. Theorem.** (Interpolation Theorem for Sets of Points) *As in Theorem 3.1 where $S = \{A, B\}$, $S^A = \{A, I\}$, $S^B = \{B, I\}$ and regular operators are defined in the weaker sense indicated in Remark 3.6.*

**Proof.** The proof is essentially a simplified version of the proof of Theorem 3.1 where we forget that structured sets are bipartite. ∎

# 4 Interpretations into logical languages

In this section we will consider Gentzen-like logical systems (such as $LK$, $LJ$, etc.) where all rules (and therefore proofs) satisfy the well-known *subformula property*

416      A. Carbone

**3.5. Remark.** While in the proof of Theorem 3.1 we need operators derived by adding only copies of the empty set, it is worthwhile to notice that in the proof of Theorem 3.4 operators are derived by using auxiliary sets which are in general non-trivial, i.e. they are $I$-sets equipped with a certain structure. This induces transformations based on derived operators which have as a value a set with a structure somehow 'dependent' on the structure of these auxiliary sets. Clearly empty sets do not impose any structural requirement. In this sense the reader can be lead to the impression that in the proof of Theorem 3.4 we use a substantially different construction than in the proof of Theorem 3.1. Indeed this is not the case since one can show that any structured set so obtained can be $\mathcal{O}$-defined from non-restricted structured sets. The basic point (and this is a general fact) is that for any $\mathcal{O}$-definable structured set $\langle S^1 | S^2 \rangle$ (obtained with the application of $p$ steps of transformation) and any $S^1_* \supseteq S^1, S^2_* \supseteq S^2$, one can $\mathcal{O}$-define the structured set $\langle S^1_* | S^2_* \rangle$ (with $p$ steps of transformation) by adding to the trivial structured sets defining $\langle S^1 | S^2 \rangle$ those sets of points in $S^1_* \backslash S^1, S^2_* \backslash S^2$ and then drag them along the transformations as required. Once we know this, we are done. In fact, suppose we can derive a structured set $\langle S^1_0 | S^2_0 \rangle$ from an $\mathcal{O}$-definable structured set $\langle S^1 | S^2 \rangle$ by applying a unary operator of subarity $l < k$ derived from a unary operator of subarity $k$ by using $k - l$ auxiliary sets. Then, we can notice that by 'adding' to the suitable components of $\langle S^1 | S^2 \rangle$ the desired $k - l$ sets (using the general fact pointed out before), one can obtain $\langle S^1_0 | S^2_0 \rangle$ by the direct application of the unary operator of subarity $k$. (What we just described has a well-known counterpart in the context of logic with the distinction between multiplicative and additive unary rules; notice that the basic property we used corresponds to the *weakening* rule (see section 4).) The same idea can be applied to derived operators of arity $k$ and subarity $l < k$ to show that $\langle S^1_0 | S^2_0 \rangle$ can be obtained by applying the $l$-ary derived operator only to auxiliary sets which are empty.

**3.6. Remark.** If we define a structured set simply as a *collection* of sets $s_1 \ldots s_n$ (equipped with some additional structure) instead of a bipartite collection, we can show an interpolation theorem for sets $S = \{A, B\}$, where regular pairs are the same as before except that we forget the bipartite structure.

**3.7. Theorem.** (Interpolation Theorem for Sets of Points) *As in Theorem 3.1 where $S = \{A, B\}$, $S^A = \{A, I\}$, $S^B = \{B, I\}$ and regular operators are defined in the weaker sense indicated in Remark 3.6.*

**Proof.** The proof is essentially a simplified version of the proof of Theorem 3.1 where we forget that structured sets are bipartite. ∎

## 4   Interpretations into logical languages

In this section we will consider Gentzen-like logical systems (such as $LK$, $LJ$, etc.) where all rules (and therefore proofs) satisfy the well-known *subformula property* (roughly speaking this means that any formula appearing in the proof, will appear as subformula of the end-sequent; for a precise definition see [10].) [6] To obtain the Craig Interpolation theorems for such systems, we need simply to interpret atomic occurrences as points in a set, formulas as sets in a structured space (where the structure is given by the tree-like form of logical formulas), sequents as structured sets and rules of inference for the logical system as operators on structured spaces. This is the rough idea of what we will present, but let us now define more precisely what is a logical language and a sequent calculus.

A *logical language* is a language containing logical symbols such as $\wedge, \vee, \neg, \supset, \exists, \forall, \ldots$, variables, constants, (first or/and higher order) function symbols, (first or/and higher order) predicate symbols, special predicate constants $\perp, \top, \ldots$, possibly metamathematical symbols like metavariables, term variables, predicate (formula) variables. *Formulas* are defined as usual; we will denote them with capital letters $A, B, C, \ldots, A_1, A_2, \ldots$. A *sequent* is a line of the form

$$A_1, \ldots, A_k \to B_1, \ldots, B_l$$

where the $A_i$'s and $B_j$'s are formulas; its intended meaning is $\bigwedge_i A_i \supset \bigvee_j B_j$. We permit $k$ and $l$ to be zero. A sequence of formulas separated by commas is a *cedent*; in the sequent above, $A_1, \ldots, A_k$ is the *antecedent* and $B_1, \ldots, B_l$ is the *succedent*. We refer to antecedents and succedents in a sequent using capital letters of the Greek alphabet. For instance $\Gamma \to \Delta$ denotes a sequent. In the following we will intend a sequence $A_1, \ldots, A_k$ to be a *multiset* of formulas, i.e. finite (possibly empty) set of formulas, in which repetitions of some formulas are admitted; the order of formulas in a multiset is not essential but for every member of the multiset the number of its occurrences is important.

A *sequent calculus* is defined by prescribing a set of sequents called *axioms* and a set of ($k$-ary) *rules of inference* of the form

$$\frac{S_1 \quad \ldots \quad S_k}{S_0}$$

where $S_0, S_1, \ldots, S_k$ are sequents and $k \geq 1$ ($S_1, \ldots, S_k$ are the *antecedents* of the rule and $S_0$ is its *consequent*.) A proof in a sequent calculus is intended to be a *tree* of sequents; each sequent must either be an axiom (in this case the sequent is labeling a leaf of the tree) or be derived by one of the rules of inference (the sequent is a label for an internal node of the tree). Every occurrence of a sequent in a proof other than the end-sequent is used exactly once as the premise of a rule of inference.

For the purpose of this section we do not need to enter the details of the formalization for the classical system $LK$ (its axioms and rules are anyway unmysterious). For a precise formulation we refer the reader to [10]. We are simply interested to

---

[6] The subformula property may seem a little odd to the nonlogician. In fact it implies that modus ponens for instance is not allowed to be used in a proof. It is worthwhile to mention here that the Gentzen Cut-Elimination Theorem (see [10]) permits one to transform proofs of the more familiar variety into proofs which satisfy the subformula property.

point out here that *logical axioms* in $LK$ are sequents of the form $A, \Gamma \to \Delta, A$ where we intend the two occurrences of $A$ on the right and left hand side of the sequent arrow $\to$ to be *logically linked*. Similarly, rules of inference in $LK$ always define precise logical relations between antecedent(s) and consequent. To illustrate the idea let us consider the $\wedge$:*right* rule in $LK$

$$\frac{\Gamma_1 \to \Delta_1, A \quad \Gamma_2 \to \Delta_2, B}{\Gamma_1, \Gamma_2 \to \Delta_1, \Delta_2, A \wedge B}$$

where $\Gamma_1, \Gamma_2, \Delta_1, \Delta_2$ are collections of formulas called *side* formulas of the rule, $A, B$ in the antecedents are the *auxiliary* formulas and $A \wedge B$ in the consequent is the *main* formula. All occurrences of formulas in the upper sequents appear in the lower sequent of the rule. In particular this is the case for the auxiliary formulas $A$ and $B$ occurring as subformulas of the main formula $A \wedge B$ of the rule. If we give a closer look to the rule [7] we notice that each occurrence of *atomic* formulas in the upper sequent of the rule has a copy on the lower sequent, and a *logical link* between the two occurrences is defined by the intended meaning of the rule. [8] Let us illustrate this point by tracing the logical links between atomic occurrences in a very simple proof in $LK$

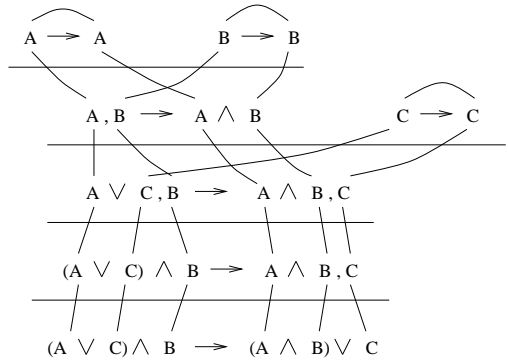

We are now ready to define the logical interpretation of the set-theoretic language to deduce from it the Craig Interpolation Theorem for $LK$.

Atomic formulas can be thought of as points in a set, and formulas as sets of points labelling the leaves of a tree-structure. The internal nodes of the tree will be labelled by the logical connectives of the formula.

---

[7] The logical relations between occurrences of formulas in a proof have been analysed in [5] and [2] where graph-theoretical notions tracing the flow of occurrences in a proof have been defined. See also [3].

[8] This fact holds not only for rules in $LK$ but for any rule of Gentzen-like systems satisfying the subformula property. Notice that here we think of a version of classical logic which does not contain the so-called *cut rule*, i.e. the rule with antecedents $\Gamma_1 \to \Delta_1, A$ and $A, \Gamma_2 \to \Delta_2$, and consequent $\Gamma_1, \Gamma_2 \to \Delta_1, \Delta_2$ (roughly speaking, a generalized form of modus ponens.) For this rule one can observe that there is no trace in $\Gamma_1, \Gamma_2 \to \Delta_1, \Delta_2$ of $A$ and in particular, of the logical link between the two occurences of the formula $A$ in the rule.

Sequents are essentially structured sets where formulas on the left hand side of the sequent arrow $\rightarrow$ belong to the first component of the partition set and those on the right hand side belong to the second component.

Axioms in $LK$ (which are of the form $A, \Gamma \rightarrow \Delta, A$) are bipartite sets containing two distinguished sets of points (corresponding to the $A$ occurrences) having obviously the same tree-structure and embedded by the identity map one into the other.

Rules in $LK$ are either unary or binary operators on bipartite sets which will take a (pair of) set(s) having a tree-structure and will form a new tree by adding a new parent to the (pair of) root(s) (if a set is empty, i.e. there is no tree structure associated to it, then we assume a point to be its associated tree.) The embeddings associated to each rule are naturally defined by those *logical links* between occurrences of atomic formulas we described above (this is illustrated in the picture above tracing pointwise the flow of the embeddings.) [9] Notice that rules are naturally interpreted as operators by the way conditions 1 and 2 in the definition of operator have been given, which ensure arguments (the auxiliary formulas of a rule) to be embedded into a value (the main formula of a rule) while the remaining sets occurring in the structured sets considered by the operator (the side formulas in the antecedents of a rule) will remain untouched.

The set of rules of $LK$ is closed under regularity conditions. This is easy to check since $LK$ contains a pair of rules for negations which move a formula from one side of the sequent arrow to the other. By combining properly conjunctions, disjunctions and negations one is able to define in $LK$ regular pairs for any given arity and choice of components for arguments and value. Notice that regular rules (i.e. rules whose interpretation is a regular pair of operators) will be either rules in $LK$ or obtained by combining rules of $LK$. Notice also that only unary and binary regular rules are needed to show that $LK$ is closed under regularity.

For what it concerns operators derived from regular rules, let us observe that unary operators of subarity 1 derived from a unary operator of subarity 2 make sense in $LK$ only when the $\vee$:*right* or $\wedge$:*left* rules are formalized with antecedents of the form $\Gamma \rightarrow \Delta, A, B$ or $A, B, \Gamma \rightarrow \Delta$ (the so-called *multiplicative* form of the rule.) In this case the derived unary operators of subarity 1 will be defined through $\vee$:*right* or $\wedge$:*left* rules with antecedents of the form $A, \Gamma \rightarrow \Delta$ or $\Gamma \rightarrow \Delta, A$ (the so-called *additive* form of the rule.) Notice that the additive form can be obtained from the multiplicative by adding the formula $B$ as extra assumption for the antecedent. This operation preserves the validity of a sequent; it is known in logic as *weakening* and it belongs to the set of rules of $LK$.

Let us also notice that a unary operator derived from a binary operator corresponds in $LK$ to a special form of inference, where one of the antecedents is an axiom either of the form $\rightarrow \top$ or $\bot \rightarrow$.

---

[9]Because of the usual requirement on eigenvariables for quantifier rules, one should obviously impose some requirement to the applicability of the operators. A detailed description leads to a technical development we will not consider here. The reader can refer to [4] where similar ideas are developed in the logical language of schematic systems.

From these two observations on derived operators, it follows that rules derived from any regular rule can be effectively constructed in $LK$. Therefore, Craig's Theorem holds for $LK$.

To conclude with the interpretation for $LK$, we should notice that in case 1 (2) of the basic step defining $\mathcal{S}_0$ in the proof of Theorem 3.1 we need just to interpret the empty set as the predicate constant $\bot$ ($\top$) if $r > v$, and $\top$ ($\bot$) if $r = v$. Then, to obtain a proof of $A \to I$ and $I \to B$ where $I$ is expressed in the language common to $A$ and $B$, we replace $\bot$ ($\top$) by any refutable (provable) formula of $LK$ whose symbols belong both to $A$ and $B$. Because of the pair of rules for negation in $LK$, notice that if the sequents $A \to I$ and $I \to B$ ($A, I \to$ and $\to B, I$) are provable then $A, \neg I \to$ and $\to B, \neg I$ ($A \to \neg I$ and $\neg I \to B$) are also provable sequents. Through this interpretation, the reader can essentially reformulate Maehara's proof for $LK$ (see p.33ff in [10].)

Let us now consider the system of intuitionistic logic $LJ$. (This is a system where we cannot prove $A \vee B$ without actually proving either $A$ or $B$; similarly we cannot prove $\exists x. A(x)$ without having already proved $A(t)$ for some term $t$. These properties fail to hold in $LK$.) The interpretation of the operators to show Craig's Theorem for intuitionistic logic $LJ$ is straightforward if we formulate $LJ$ as a system containing axioms of the form $\Gamma, C \to C$ and rules of the form

$\neg : left \quad \dfrac{\Gamma \to A}{\neg A, \Gamma \to} \qquad \neg : right \quad \dfrac{A, \Gamma \to}{\Gamma \to \neg A}$

$\wedge : right \quad \dfrac{\Gamma_1 \to A \quad \Gamma_2 \to B}{\Gamma_{1,2} \to A \wedge B}$

$\wedge : left \quad \dfrac{A, \Gamma \to C^1, \ldots, C^k}{A \wedge B, \Gamma \to C^1, \ldots, C^k} \qquad \dfrac{A, \Gamma \to C^1, \ldots, C^k}{B \wedge A, \Gamma \to C^1, \ldots, C^k}$

$\vee : left \quad \dfrac{A, \Gamma_1 \to C^1, \ldots, C^k \quad B, \Gamma_2 \to C^1, \ldots, C^l}{A \vee B, \Gamma_{1,2} \to C^1, \ldots, C^{k+l}}$

$\vee : right \quad \dfrac{\Gamma \to A}{\Gamma \to A \vee B} \qquad \dfrac{\Gamma \to A}{\Gamma \to B \vee A}$

$\forall : right \quad \dfrac{\Gamma \to A(a)}{\Gamma \to \forall x. A(x)} \qquad \forall : left \quad \dfrac{A(t), \Gamma \to C^1, \ldots, C^k}{\forall x. A(x), \Gamma \to C^1, \ldots, C^k}$

$\exists : right \quad \dfrac{\Gamma \to A(t)}{\Gamma \to \exists x. A(x)} \qquad \exists : left \quad \dfrac{A(a), \Gamma \to C^1, \ldots, C^k}{\exists x. A(x), \Gamma \to C^1, \ldots, C^k}$

$Contraction \quad \dfrac{\Gamma \to C^1, \ldots, C^s}{\Gamma \to C^1, \ldots, C^{s-1}} \qquad \dfrac{A, A, \Gamma \to C^1, \ldots, C^k}{A, \Gamma \to C^1, \ldots, C^k}$

where the $C^i$'s on the right hand side of the sequents, denote distinguished occurrences (with $k, l \geq 0, s \geq 2$) of the formula $C$. The variable $a$ (called the *eigenvariable* of the ∃:*left*/∀:*right* rule) does not occur free in $\Gamma, C$, and $t$ is an arbitrary term.

By interpreting sequents, formulas, atomic formulas and rules in the same way as for $LK$, we obtain that axioms and sequents appearing in rules of $LJ$ are restricted structured sets (see Remark 3.3.) The discussion about derived rules for $LK$ applies to $LJ$ as well. One should notice moreover that the *contraction* rules in the formalization of $LJ$ correspond to the operators of contraction in the extended closure condition required by Theorem 3.4. The Interpolation Theorem for $LJ$ follows readily from Theorem 3.4.

**4.1. Remark.** In $LK$ and $LJ$ each pair of rules introducing a connective to the right or to the left of the sequent arrow, forms a regular pair. A set of operators closed under regularity might contain operators which do not belong to any regular pair though. This means that Craig's Theorem can be shown for logical systems of rules which are closed under the regularity condition, but that contain also rules not belonging to any regular pair. Most of the logics usually studied in the literature are defined exactly by pairs of regular rules though. This fact was also observed in [4] and used to deduce Craig's Theorem for certain fragments of linear logic and fragments of classical logic.

**4.2. Remark.** Theorem 3.1 and its variants (i.e. Theorem 3.4 and Theorem 3.7) when applied to logic point out that a logical system does not need to satisfy the *subformula property* for its rules in order to enjoy Interpolation, but simply some form of *embedding* of the antecedents of a rule into its consequent. In particular, our definition of embedding can be weakened so to interpret logical systems formalizing proofs with analytic cuts (i.e. cuts on formulas appearing as subformulas of the endsequent of the proof) which are well known to enjoy Interpolation (the reader can refer to [3] for results in this direction.)

**4.3. Remark.** Since the definition of an operator is based on embeddings of *arbitrary* nature, we have that any pair of sets (formulas) $A$ and $B$ definable (built) by our operators (rules) is not related by a map but by a pseudomap. Even though $LK$ and $LJ$ are in a sense very simple systems (their rules are either binary or unary, and belong to regular pairs), we observe the same phenomena also for their interpretation. This complexity of $LK$ and $LJ$ is due to the *contraction* rule in their formalization.

# 5  Discussion

As a result of our analysis we have now a number of conditions sufficient for a system of combinatorial nature to enjoy Interpolation. Its objects might be graphs just as well as formulas or surfaces. Let us step back for a moment and take a wider look at our combinatorial model. Let $\Omega$ denote the set of all spaces $\mathcal{S}$ defined with respect to some universe, with some choice of structure for sets of points and some prescribed class of embeddings. There is a distinguished subset of $\Omega$, the set $\Omega_0$ of trivial spaces. A collection of operators $\mathcal{O}$ induces a collection of pseudomaps on $\Omega$. We can think of $\Omega$ together with this collection of pseudomaps as defining a sort of dynamical system. Through this view, the object of our analysis turns out to be the orbit of $\Omega_0$ under the collection of transformations. This is not a dynamical system in the more typical sense, because our transformations are pseudomappings instead of mappings. In general if pseudomappings were allowed in dynamical systems, one should not expect to have as much interesting structure as usual. On the other hand, the use of pseudomappings in logic is balanced by the fact that one defines them through a finite number of operators which are fairly simple and encode sufficient additional structure (i.e. the logical structure of the formulas) to still have interesting combinatorics.

Witness of the intriguing combinatorics encoded in logical proofs is the link between the complexity of an interpolant and more classical questions in complexity theory. For classical logic with equality, it has been proved by Meyer ([8]) that there is no general recursive bound on the length of the smallest interpolant. It is an open question whether or not the size (i.e. the number of symbols) of the smallest interpolant in the propositional case (roughly, classical logic free of quantifiers) can be polynomially bounded by the size of the tautology. A positive answer to this open question would have as a consequence in complexity theory that $NP \cap co-NP \subset P/poly$ (for recent work on this direction see [3], [6], [7].) To see this consider a language $L \subset \{0,1\}^*$ and $L \in NP \cap co-NP$ (therefore the complement $\bar{L}$ of $L$ is in $NP$). It is well-known that there are sequences of propositional formulas $A_n(p_1, \ldots, p_n, q_1, \ldots, q_m), B_n(p_1, \ldots, p_n, r_1, \ldots, r_s)$ (where $p_1, \ldots, p_n$, $q_1, \ldots, q_m$, $r_1, \ldots, r_s$ are the only propositional variables occurring in $A_n, B_n$) such that their size is $|A_n|, |B_n| \leq n^{\mathcal{O}(1)}$ and

$$L = \{(p_1, \ldots, p_n) \in \{0,1\}^n \mid \exists y_1 \ldots y_m . A_n(p_1, \ldots, p_n, y_1, \ldots, y_m) \; holds\}$$

$$\bar{L} = \{(p_1, \ldots, p_n) \in \{0,1\}^n \mid \exists z_1 \ldots z_s . B_n(p_1, \ldots, p_n, z_1, \ldots, z_s) \; holds\}$$

Clearly $A_n \rightarrow \neg B_n$ is a tautology (for all $n$). By hypothesis there is an interpolant $I_n(p_1, \ldots, p_n)$ for $A_n \rightarrow \neg B_n$ such that $|I_n| \leq f(|A_n| + |B_n|)$ (for all $n$) where $f$ is some function computable in polynomial time. Therefore $L \in P/poly$.

Since the size of an interpolant $I_n$ constructed in Theorem 3.1 depends linearly on the size of a cut-free proof of $A_n \rightarrow \neg B_n$ (for all $n$), one can wonder where the complexity problem hides for the provability of tautologies. There are families of tautologies in propositional logic which can be proved using a polynomial number of symbols with the help of the cut rule but only with an exponential number of symbols if the cut-rule is not allowed. An example is the family of propositional tautologies expressing the pigeon hole principle ([1].) Families of tautologies which are known to behave this way, have usually linear interpolants. In [9] one finds a quadratic lower bound. As we mentioned in footnote 6 one can always transform a proof (of polynomial size) containing cuts into a proof without cuts (with exponential size in the worse case.) The combinatorial meaning of the procedure of cut-elimination is by no means clear.

*Institut für Algebra und Diskrete Mathematik*
*Technische Universität Wien*
*Wiedner Hauptstrasse 8-10/118*
*A-1040 Wien*
*Austria*

e-mail: ale@logic.tuwien.ac.at